\documentclass[a4paper]{article}
\usepackage{graphicx}
\usepackage{amsmath}
\usepackage{amsthm}
\usepackage{amssymb}
\usepackage{cite}

\newcommand{\supp}{\mathop{\mathrm{supp}}\nolimits}

\newtheorem{theorem}{Theorem}
\newtheorem{lemma}{Lemma}
\newtheorem{corollary}{Corollary}
\newtheorem{fact}{Fact}

\newtheorem{remark}{Remark}
\newtheorem{example}{Example}
\newtheorem*{acknowledgements}{Acknowledgements}

\title{Statistical properties of a generalized threshold network model}
\author{
Yusuke Ide,\ Norio Konno
\thanks{Postal address: 79-5, Tokiwadai, Hodogaya-ku, Yokohama 240-8501, Japan}
\\
\textit{Department of Applied Mathematics,}
\\
\textit{Yokohama National University}
\\
\and
Naoki Masuda 
\thanks{Postal address: 7-3-1, Hongo, Bunkyo-ku, Tokyo 113-8656, Japan}
\\
\textit{Graduate School of Information Science and Technology,}
\\
\textit{The University of Tokyo}
}
\date{}
\begin{document}
\maketitle
\begin{abstract}
The threshold network model is a type of finite random graphs. In this
paper, we introduce a generalized threshold network model. A pair of 
vertices with random weights is connected by an edge when real-valued 
functions of the pair of weights belong to given Borel sets. 
We extend several known limit theorems for the number of prescribed 
subgraphs to show that the strong law of large numbers can be uniform convergence. 
We also prove two limit theorems for the local and global clustering coefficients.
\end{abstract}
\section{Introduction}
Complex networks have been an attractive research topic for a
decade. Particularly, many real-world graphs are characterized by
the small diameter, high clustering (abundance of connected 
triangles), and fat-tail degree distributions. 
Degree distributions often follow the truncated 
power law, which is called the scale-free property of
networks
\cite{Albert02,NewmanSIAM,blmch}. Both growing and static network models
are capable of generating scale-free networks.

Here we are concerned with asymptotic properties of a class of static
network models called the threshold network model, which is generated on
$n$ vertices labeled $1,\dots ,n$ with independent and
identically distributed (i.i.d.)\ random weights $X_{1},\ldots
,X_{n}$. We connect a pair of vertices $i$ and $j$ with $i\neq j$ by
an edge when $X_{i}+X_{j}>\theta $ for a given threshold $\theta$.
The threshold network model is a subclass of so called hidden variable
models and its mean behavior \cite{bps,ccdm,hss,mmk04,mmk05,scb,so} and limit
theorems \cite{kmrs05,nr03} have been analyzed.

To define a generalization of the threshold network model,
let $\mathbb{R}^{d}$ be the $d$-dimensional Euclidean space. 
We prepare an i.i.d.\ sequence of $\mathbb{R}^{d}$-valued random variables 
$X_{1},\dots ,X_{n}$ with a common distribution function $F$. 
We associate the random variable $X_{i}$, which we call
weight function,  with vertex $i$.
Now we introduce Borel measurable functions 
$f_{c}^{ l^{\prime } }: (\mathbb{R}^{d})^{2}\to \mathbb{R}$ 
with $f_{c}^{ l^{\prime } }(x,y)=f_{c}^{ l^{\prime } }(y,x)$
for all $l^{\prime }\in \{1,\dots ,l\}$. 
Let $\mathcal{B}(\mathbb{R})$ be the Borel $\sigma $-field of
$\mathbb{R}$.  For a given finite collection of Borel measurable sets
$\mathcal{C}=\{B_{1},\dots ,B_{l}\}$ with $B_{ l^{\prime } }\in
\mathcal{B}(\mathbb{R})$,
we connect vertices
$i$ and $j$ ($i\neq j$) if $f_{c}^{ l^{\prime } }(X_{i},X_{j})\in
B_{ l^{\prime } }$ for all ${ l^{\prime } }\in \{1,\dots ,l\}$. In
other words, we form an edge $\langle i,j\rangle $ if
$
\prod _{l^{\prime }=1}^{l}I_{B_{ l^{\prime } }}
\left(f_{c}^{ l^{\prime } }(X_{i},X_{j})\right)=1 
$
for $i\neq j$,
where $I_{A}(x)$ denotes the indicator function, i.e., 
$I_{A}(x)=1$ for $x\in A$ and $I_{A}(x)=0$ otherwise. 
Thus we obtain a random graph $G_{\mathcal{C}}(X_{1},\dots ,X_{n})$. 
If there exist two collections of Borel sets $\mathcal{C}=\{B_{1},\dots ,B_{l}\}$ and 
$\mathcal{C}^{\prime }=\{B^{\prime }_{1},\dots ,B^{\prime }_{l}\}$ 
with $B_{ l^{\prime } }\subset B^{\prime }_{ l^{\prime } }$ for all $l^{\prime }\in \{1,\dots ,l\}$, then 
$\mathbb{P}\{\langle i,j\rangle \in G_{\mathcal{C}}(X_{1},\dots ,X_{n})\}\leq 
\mathbb{P}\{\langle i,j\rangle \in G_{\mathcal{C}^{\prime }}(X_{1},\dots ,X_{n})\}$ 
holds by simple coupling.

This random graph generalizes the threshold network
model studied in \cite{bps,ccdm,hss,kmrs05,mmk04,mk06,nr03}. 
By choosing $l=1$, $B_1 = (\theta ,\infty )$ for some 
$\theta \in \mathbb{R}$, $f_{c}^1(x,y)=x+y$,
we reproduce the model in \cite{bps,ccdm,hss,kmrs05,mmk04}. 
In the context of social networks, 
a model with $l=2$,
$B_{1}=(\theta ,\infty )$, $B_{2}=(-\infty ,c]$ 
($\theta ,c\in \mathbb{R}$),
$f_{c}^{1}(x,y)=x+y$, and $f_{c}^{2}(x,y)=|x-y|$ (or
$f_{c}^{2}(x,y)=|x-y|/(x+y)$) was proposed \cite{mk06}.
General limit theorems are shown in \cite{nr03} when
$l=1$, $B_1 = (\theta ,\infty )$, and 
$f_{c}^{1}(x,y)=|x-y|$.

In Sec.\ 2, 
we state several general limit theorems for the number of
prescribed subgraphs.  By using $U$-statistics, the strong law of
large numbers, the central limit theorem, and the law of the iterated
logarithm are stated for global properties of the model.  We
also state a limit distribution for a local property. 
These are generalizations of Theorems 1, 4, and 5 of \cite{kmrs05} and 
Theorems 1, 2, and 3 of \cite{nr03}. 
In Sec.\ 3, 
we show that the strong law of large numbers for the number of prescribed
subgraphs is uniform convergence on so-called the VC class of Borel sets,
generalizing Theorem 1a of \cite{nr03}. 
In Sec.\ 4, 
we show limit theorems for the clustering coefficient, which quantifies
the abundance of connected triangles in a graph in a specific ways. 
Particularly, we show the strong law of large numbers for the local clustering coefficient (Theorem \ref{thmloccluster}) 
and the global clustering coefficient (Theorem \ref{thmglobcluster}). 
Theorems \ref{thmloccluster} and \ref{thmglobcluster} are main results of this paper. 
In Sec.\ 5, 
we present several examples of limit degree distributions.
\section{General Results}
In this section, we show limit theorems
for the number of prescribed subgraphs. 
Let us begin with notations \cite{nr03}. 
For $m\in \{2,\dots ,n\}$, we consider a graph $H=(V_{H},E_{H})$ on
the ordered set of $m$ vertices $V_{H}=(v_{1},\dots ,v_{m})$ and the
edge set $E_{H}$. For another graph $H^{\prime }=(V_{H^{\prime
}},E_{H^{\prime }})$ on $m$ vertices, we say $H^{\prime }\thicksim H$
if $V_{H^{\prime }}=V_{H}$ and $E_{H^{\prime }}=E_{H}$ for some
reordering of vertices.  Thus $\mathcal{A}_{H}^{\thicksim
}=\{H^{\prime }:H^{\prime }\thicksim H\}$ is the
set of all graphs isomorphic to $H$.  Let us define
$\mathcal{A}_{m}=\bigcup _{i}\mathcal{A}_{H_{i}}^{\thicksim }$, where
$H_i$ is an arbitrarily chosen graph on $m$ vertices.
The collection of all triangles 
and graphs on three vertices that consist of two connected vertices
and an isolated vertex 
is an example of $\mathcal{A}_{3}$. 
The collection of cliques on $m$ vertices 
and the graphs on $m$ vertices with $m$ isolated vertices 
is an example of $\mathcal{A}_{m}$. 
The definition of $\mathcal{A}_{m}$ guarantees the
symmetrical property of the kernel function 
$h_{\mathcal{A}_{m}}:(\mathbb{R}^{d})^{m}\to \mathbb{R}$ given by 
\begin{eqnarray}\label{defkernel}
h_{\mathcal{A}_{m}}(x_{1},\dots ,x_{m})
=
I_{\mathcal{A}_{m}}\left(G_{\mathcal{C}}(x_{1},\dots ,x_{m})\right),
\end{eqnarray}
where $G_{\mathcal{C}}(x_{1},\dots ,x_{m})$ denotes a realization of 
the random graph $G_{\mathcal{C}}(X_{1},\dots ,X_{m})$. 
Then we define 
\begin{eqnarray*}
\Tilde{U}_{n}(\mathcal{C},\mathcal{A}_{m})=
\sum _{1\leq i_{1}<\dots <i_{m}\leq n}
h_{\mathcal{A}_{m}}(X_{i_{1}},\dots ,X_{i_{m}}),
\end{eqnarray*}
i.e., the number of 
subgraphs belonging to the collection $\mathcal{A}_{m}$
in the random graph $G_{\mathcal{C}}(X_{1},\dots ,X_{n})$.
We also define 
\begin{eqnarray*}
U_{n}(\mathcal{C},\mathcal{A}_{m};i)=
\sum _{
\begin{subarray}{c}
1\leq i_{2}<\dots <i_{m}\leq n\\
i_{2},\dots ,i_{m}\neq i
\end{subarray}
}
h_{\mathcal{A}_{m}}(X_{i},X_{i_{2}},\dots ,X_{i_{m}}),
\end{eqnarray*}
i.e., \ the number of subgraphs that include vertex $i$ 
and belong to $\mathcal{A}_{m}$
in the random graph $G_{\mathcal{C}}(X_{1},\dots ,X_{n})$. 

Note that $U_{n}(\mathcal{C},\mathcal{A}_{m};i), 1\leq i\leq n$ are
identical in distribution and the following relation holds:
\begin{eqnarray*}
\frac {\sum _{i=1}^{n}U_{n}(\mathcal{C},\mathcal{A}_{m};i)/\binom{n-1}{m-1}}{n}
=\frac{m\Tilde{U}_{n}(\mathcal{C},\mathcal{A}_{m})}{n\binom{n-1}{m-1}}
=\frac{\Tilde{U}_{n}(\mathcal{C},\mathcal{A}_{m})}{\binom {n}{m}}.
\end{eqnarray*}
This implies that the global property 
$\Tilde{U}_{n}(\mathcal{C},\mathcal{A}_{m})/\binom {n}{m}$ 
is the arithmetic mean of the local properties 
$U_{n}(\mathcal{C},\mathcal{A}_{m};1)/\binom {n-1}{m-1},\dots ,
U_{n}(\mathcal{C},\mathcal{A}_{m};n)/\binom {n-1}{m-1}$. 

We define 
\begin{eqnarray*}
F(\mathcal{C},\mathcal{A}_{m})
\! \! \! \! &=&\! \! \! \!
\mathbb{E}[h_{\mathcal{A}_{m}}(X_{1},\dots ,X_{m})],\\
\zeta (\mathcal{C},\mathcal{A}_{m})
\! \! \! \! &=&\! \! \! \! 
Var(\mathbb{E}[h_{\mathcal{A}_{m}}(X_{1},\dots ,X_{m})|X_{1}]),
\end{eqnarray*}
and assume $\zeta (\mathcal{C},\mathcal{A}_{m})>0$. 
Since $\Tilde{U}_{n}(\mathcal{C},\mathcal{A}_{m})/\binom{n}{m}$ is a $U$-statistic \cite{s} obtained from 
the symmetric kernel $h_{\mathcal{A}_{m}}$, 
the strong law of large numbers (SLLN), the central limit theorem (CLT) and the law of the iterated logarithm (LIL) 
are derived from general results for the $U$-statistics, namely,
Theorem A (SLLN) and Theorem B (LIL) in Section 5.4, and Theorem A
(CLT) in Section 5.5 of \cite{s}: 
\begin{fact}[SLLN for global property]\label{thmglobu}As $n\to \infty $, 
\begin{eqnarray*}
\frac {\Tilde{U}_{n}(\mathcal{C},\mathcal{A}_{m})}{\binom {n}{m}}
\to 
F(\mathcal{C},\mathcal{A}_{m}),\quad \text{almost surely}. 
\end{eqnarray*}
\end{fact}
\begin{fact}[CLT for global property]As $n\to \infty $,
\begin{eqnarray*}
\displaystyle \sqrt {\frac{n}{m^{2}\zeta (\mathcal{C},\mathcal{A}_{m})}}
\Biggl[\frac {\Tilde{U}_{n}(\mathcal{C},\mathcal{A}_{m})}{\binom {n}{m}}-F(\mathcal{C},\mathcal{A}_{m})\Biggr]
\Longrightarrow \mathcal{Z}, 
\end{eqnarray*}
where $\Longrightarrow $ stands for convergence in distribution and 
$\mathcal{Z}$ is a standard normal random variable. 
\end{fact}
\begin{fact}[LIL for global property]
\begin{eqnarray*}
\limsup _{n\to \infty }\sqrt {\frac{n(\log \log n)^{-1}}{2m^{2}\zeta (\mathcal{C},\mathcal{A}_{m})}}
\Biggl|\frac {\Tilde{U}_{n}(\mathcal{C},\mathcal{A}_{m})}{\binom {n}{m}}-F(\mathcal{C},\mathcal{A}_{m})\Biggr|=1,
\quad \text{almost surely}.
\end{eqnarray*}
\end{fact}
There are the direct generalization of Theorem 4 of \cite{kmrs05} and 
Theorems 1, 2, and 3 of \cite{nr03} 
to the present model. 
\begin{remark}
Generally, when the $f^{ l^{\prime } }_{c}$ is asymmetric (e.g. directed graph), 
the number of graphs isomorphic to a graph $H$ on $m$ vertices is $\binom{n}{m}\cdot m!$. 
The limit theorems above are valid by replacing the normalizing factor $\binom{n}{m}$ 
with $\binom{n}{m}\cdot m!$. 
\end{remark}
By generalizing Theorems 1 and 5 of \cite{kmrs05}, we obtain the following 
asymptotic behavior:
\begin{fact}\label{thmlocu}As $n\to \infty $,
\begin{eqnarray*}
\frac {U_{n}(\mathcal{C},\mathcal{A}_{m};1)}{\binom {n-1}{m-1}}
\Longrightarrow U(\mathcal{C},\mathcal{A}_{m}),
\end{eqnarray*}
where 
\begin{eqnarray*}
U(\mathcal{C},\mathcal{A}_{m})=\int _{\mathbb{R}^{d}}\! \!\dotsi \! \! \int _{\mathbb{R}^{d}}
h_{\mathcal{A}_{m}}(X_{1},x_{2},\dots ,x_{m})F(dx_{2})\cdots F(dx_{m}). 
\end{eqnarray*}
\end{fact}
\begin{remark}
Almost sure convergence theorem for 
$U_{n}(\mathcal{C},\mathcal{A}_{m};1)/\binom {n-1}{m-1}$
in the sense of Theorem \ref{thmloccluster} in Sec.\ 4.1 can be proved by a simple modification 
of the proof of Theorems 1 and 5 of \cite{kmrs05}. 
\end{remark}
\section{Uniform Property}
The Vapnik-Chervonenkis approach is well known in the context of the statistical learning theory. 
Particularly, it is useful in showing uniform convergence for limit theorems \cite{d,p}. 
In this section, we show that SLLN for global property (Fact \ref{thmglobu}) 
is uniform convergence on the VC class of the Borel sets,
which extends the special case treated in \cite{nr03}.

Let $M$ be a set and $\mathcal{D}$ be a family of subsets of $M$. 
For $A\subset M$ let $\Delta ^{\mathcal{D}}(A)=\sharp (A\cap \mathcal{D})$, 
where $\sharp (A\cap \mathcal{D})$ denotes the number of sets in 
$A\cap \mathcal{D}=\left\{A\cap D : D\in \mathcal{D}\right\}$. 
Let 
$m_{\mathcal{D}}(n)=\max _{A\subset M}\left\{\Delta ^{\mathcal{D}}(A) : |A|=n\right\}$
for $n=0,1,2,\ldots $, 
where $|A|$ denotes the number of elements in $A$, 
or if $|M|<n$ let $m_{\mathcal{D}}(n)=m_{\mathcal{D}}\left(|M|\right)$. 
We define an indicator of the family $\mathcal{D}$: 
\begin{eqnarray*}
S(\mathcal{D})
=
\begin{cases}
\sup \left\{n: m_{\mathcal{D}}(n)=2^{n}\right\} 
& \text{if $\mathcal{D}$ is non-empty},\\
-1
& \text{if $\mathcal{D}$ is empty}.
\end{cases}
\end{eqnarray*}
The family $\mathcal{D}$ is called 
a Vapnik-Chervonenkis (VC) class of sets 
if $S(\mathcal{D})<+\infty $. 
For example, the collection of half intervals $\mathcal{D}=\left\{(-\infty ,x] : x\in \mathbb{R}\right\}$ 
is a VC class on $\mathbb{R}$ with $S(\mathcal{D})=1$. 
Based on Chapter 4.5 of \cite{d}, we have  
\begin{corollary}\label{corvc}
For any $\mathcal{D}\subset 2^{M}$ and $\mathcal{D}^{\prime }\subset 2^{M}$(resp. $2^{N}$), 
if $S(\mathcal{D})<\infty $ and $S(\mathcal{D}^{\prime })<\infty $ then 
$S(\mathcal{D}\cup \mathcal{D}^{\prime })<\infty $ and $S(\mathcal{D}\cap \mathcal{D}^{\prime })<\infty $
(resp. $S(\mathcal{D}\times \mathcal{D}^{\prime })<\infty $) \\
where 
$\mathcal{D}\cup \mathcal{D}^{\prime }=\left\{D\cup D^{\prime } : 
D\in \mathcal{D},\ D^{\prime }\in \mathcal{D}^{\prime }\right\}$, 
$\mathcal{D}\cap \mathcal{D}^{\prime }=\left\{D\cap D^{\prime } : 
D\in \mathcal{D},\ D^{\prime }\in \mathcal{D}^{\prime }\right\}$
and 
$\mathcal{D}\times \mathcal{D}^{\prime }=\left\{D\times D^{\prime } : 
D\in \mathcal{D},\ D^{\prime }\in \mathcal{D}^{\prime }\right\}$. 
\end{corollary}
For a given function $h : M\rightarrow \mathbb{R}$, the subgraph of $h$ is the set 
\begin{eqnarray*}
D_{h}=\left\{(x,t)\in M\times \mathbb{R} : 0\leq t\leq h(x)\ \text{or}\ h(x)\leq t \leq 0\right\}.
\end{eqnarray*}
A class of functions $\mathcal{H}$ is a VC-subgraph class if the collection 
$\mathcal{D}_\mathcal{H}=\left\{D_{h} : h\in \mathcal{H}\right\}$ is a VC class of sets. 

For a class of real-valued measurable functions $\mathcal{H}$
on $M^{m}$ for a fixed integer $m$, Arcones and Gin\'{e} \cite{ag} proved
the following uniform SLLN for i.i.d. sequence $\left\{X_{i}\right\}_{i=1,2,\ldots }$ on $M$: 
\begin{lemma}[Corollary 3.3 of \cite{ag}]\label{lemvcag}
If $\mathcal{H}$ is a measurable VC-subgraph class of functions with \\
$\mathbb{E} \left[\sup _{h\in \mathcal{H}}|h(X_{1},\dots ,X_{m})|\right]<+\infty $, 
then 
\begin{eqnarray*}
\sup _{h\in \mathcal{H}}
\left[
\left|
\frac{1}{\binom{n}{m}}
\sum _{1\leq i_{1}<\dots <i_{m}\leq n}h(X_{i_{1}},\dots ,X_{i_{m}})
-
\mathbb{E}\left[h(X_{1},\dots ,X_{m})\right]
\right|
\right]
\to 0,\ \text{almost surely},
\end{eqnarray*}
as $n\to \infty $.
\end{lemma}

In order to use Lemma \ref{lemvcag}, 
we rewrite the kernel function (\ref{defkernel}) 
and show that the family of the kernel functions indexed 
by the VC class of Borel sets is a VC-subgraph class. 
We assume $d=l=1$. 
For a fixed integer $m\geq 2$ and
$f_{c}^1 : \mathbb{R}^{2}\to \mathbb{R}$, 
we define a function $G : \mathbb{R}^{m}\to \mathbb{R}^{\binom{m}{2}}$ by 
\begin{eqnarray*}
G(x_{1},\dots ,x_{m})=
\left(
f_{c}^1(x_{1},x_{2}),f_{c}^1(x_{1},x_{3}),\dots ,f_{c}^1(x_{1},x_{m}),
f_{c}^1(x_{2},x_{3}),\dots ,f_{c}^1(x_{m-1},x_{m})
\right).
\end{eqnarray*}
Each coordinate corresponds to a potential edge of the graph 
with a lexicographic order. 
For example, if $m=4$, the first coordinate corresponds to the 
$\langle 1,2\rangle$, 
the second coordinate to $\langle 1,3\rangle$, and 
the sixth coordinate to $\langle 3,4\rangle$. 
Note that edge $\langle i,j\rangle$
exists if and only if $f_{c}^1(x_{i},x_{j})\in B$ for a given Borel set $B$.

For a given collection
$\mathcal{A}_{m}$ of graphs on $m$ vertices, we define a set
$\mathcal{\tilde{A}}_{m}$ on $\mathbb{R}^{\binom{m}{2}}$ as follows.
For each graph $H\in \mathcal{A}_{m}$, we associate a set $\tilde{H}$
on $\mathbb{R}^{\binom{m}{2}}$. 
When a pair of vertices in
$H$ has an edge, the corresponding
coordinate of $\tilde{H}$ is occupied by $B$. Otherwise,
it is occupied by $B^{c}$, where $B^{c}$ denotes the complement of $B$.
For example, if $m=4$ and $H$ has edge set $\left\{\langle
1,2\rangle ,\langle 2,3\rangle \right\}$, then $\tilde{H}=B\times
B^{c}\times B^{c}\times B\times B^{c}\times B^{c}$.  Then we define
the set $\mathcal{\tilde{A}}_{m}=\bigcup _{H\in
\mathcal{A}_{m}}\tilde{H}$.  Note that $G(x_{1},\ldots ,x_{m})\in
\mathcal{\tilde{A}}_{m}$ is equivalent to the event that the realization
of the random graph with weights $x_{1},\ldots ,x_{m}$ is in
$\mathcal{A}_{m}$.  Finally, we obtain the rewritten
form of the kernel function:
\begin{eqnarray}\label{defkernel2}
h_{\mathcal{A}_{m}}^{B}(x_{1},\dots ,x_{m})
=
I_{\mathcal{\tilde{A}}_{m}}\left(G(x_{1},\dots ,x_{m})\right)
=
I_{G^{-1}(\mathcal{\tilde{A}}_{m})}(x_{1},\dots ,x_{m}),
\end{eqnarray}
where $G^{-1}(\mathcal{\tilde{A}}_{m})$ denotes the inverse image of
$\mathcal{\tilde{A}}_{m}$. 

Let $\mathcal{D}$ be a VC class of Borel sets on $\mathbb{R}$. 
For a fixed collection $\mathcal{A}_{m}$, 
we consider the class of kernel functions 
$\mathcal{H}_{\mathcal{D}}=\left\{h_{\mathcal{A}_{m}}^{B}(x_{1},\dots ,x_{m}) : B\in \mathcal{D}\right\}$. 
By Corollary \ref{corvc}, 
if $\mathcal{D}$ is a VC class of sets in general,
then the class of indicators 
$\left\{I_{D} : D\in \mathcal{D}\right\}$ is a VC-subgraph class. 
From Eq. (\ref{defkernel2}), if the collection 
$G^{-1}_{\mathcal{D}}(\mathcal{\tilde{A}}_{m})=\{G^{-1}(\mathcal{\tilde{A}}_{m}) : B\in \mathcal{D}\}$ 
is a VC class then $\mathcal{H}_{\mathcal{D}}$, 
which is a collection of indicator functions, 
is a VC-subgraph class on $\mathbb{R}^{m}$. 
Indeed, $G^{-1}_{\mathcal{D}}(\mathcal{\tilde{A}}_{m})$ is a VC class
based on Theorem 4.2.3 of \cite{d} and Corollary \ref{corvc} above, 
and $\mathcal{H}_{\mathcal{D}}$ is a VC-subgraph class. 
Finally, we find the uniform version of Fact \ref{thmglobu} by Lemma \ref{lemvcag}:
\begin{theorem}\label{uniform}
If $\mathcal{D}$ be a VC class of Borel sets, then for a fixed $f_{c}$ and $\mathcal{A}_{m}$, 
\begin{eqnarray*}
\sup _{B\in \mathcal{D}}
\left[
\left|
\frac{1}{\binom{n}{m}}
\sum _{1\leq i_{1}<\dots <i_{m}\leq n}h_{\mathcal{A}_{m}}^{B}(X_{i_{1}},\dots ,X_{i_{m}})
-
\mathbb{E}\left[h_{\mathcal{A}_{m}}^{B}(X_{1},\dots ,X_{m})\right]
\right|
\right]
\to 0,
\end{eqnarray*}
almost surely, as $n\to \infty $.
\end{theorem}
Theorem \ref{uniform} can be extended to general $d$ and $l$.
It generalizes
 the uniform SLLN in Theorem 1a of \cite{nr03}, which deals with
the collection of half intervals as a
VC class of sets.
\section{Clustering Coefficient}
Real-world networks are often equipped with high clustering,
that is, a large number of connected triangles.
The clustering coefficient quantifies the density of triangles in a graph
(see \cite{Albert02,NewmanSIAM} for review).
In this section, we study the limit theorems for the clustering coefficient. 
\subsection{Local Clustering Coefficient}
We assume $d=l=1$; extensions of the following results to
general $d$ and $l$ is straightforward.
We consider
a random graph $G_{B}(X_{1},\dots ,X_{n})$
for a given Borel set $B$ and $f_{c}\equiv f_{c}^{1}$.
Then we define 
\begin{eqnarray*}
D_{n}(i)
\! \! \! \! &=&\! \! \! \! 
\sum _{
\begin{subarray}{c}
1\leq j\leq n\\
j\neq i
\end{subarray}
}
h_{D}(X_{i},X_{j}),
\\
T_{n}(i)
\! \! \! \! &=&\! \! \! \! 
\sum _{
\begin{subarray}{c}
1\leq j<k\leq n,\\
j,k\neq i
\end{subarray}
}
h_{T}(X_{i},X_{j},X_{k}),
\end{eqnarray*}
where $h_{D}(x,y)=I_{B}(f_{c}(x,y))$ and 
$h_{T}(x,y,z)=I_{B}(f_{c}(x,y))\cdot I_{B}(f_{c}(y,z))\cdot I_{B}(f_{c}(z,x))$, 
i.e., $D_{n}(i)$ is the degree of vertex $i$ and $T_{n}(i)$ is the number of triangles 
including vertex $i$.
The local clustering coefficient $C_{n}(i)$ of vertex $i$ is given
by
\begin{eqnarray*}
C_{n}(i)=
\frac{ T_{n}(i) }{ \binom { D_{n}(i) }{2} }\cdot I_{\{D_{n}(i)\geq 2\}}+
w\cdot I_{\{D_{n}(i)=0,1\}}
\end{eqnarray*}
for an indeterminate $w$. 
The second term represents the
singular part for which the local clustering coefficient is not
defined in physics literature and 
applications. Here we retain this term to assess the contribution of 
vertices with degree 0 or 1. 
If it is necessary to restrict $C_{n}(i)\in [0,1]$, 
we must substitute a real value on $[0,1]$ into $w$. 
If we substitute $0$ into $w$,
the contribution of these vertices to $C_n(i)$ 
is ignored. If we substitute $1$, this contribution is
implied to be the maximum because vertices with degree
more than one satisfies $C_n(i)\le 1$.
Now we define 
\begin{eqnarray*}
V_{n}(i)=
\sum _{
\begin{subarray}{c}
1\leq j<k\leq n,\\
j,k\neq i
\end{subarray}
}
h_{V}(X_{i},X_{j},X_{k}),
\end{eqnarray*}
where $h_{V}(x,y,z)=I_{B}(f_{c}(x,y))\cdot I_{B}(f_{c}(x,z))$, 
which represents the number of vertex pairs $(j,k)$ such that both
vertex $j$ and vertex $k$ are connected to
vertex $i$.
We note the relation: On $\{D_{n}(i)\geq 2\}$ or equivalently $\{V_{n}(i)\geq 1\}$, 
\begin{eqnarray*}
\binom { D_{n}(i) }{2}=V_{n}(i),
\end{eqnarray*}
which leads to
\begin{eqnarray*}
C_{n}(i)=
\frac{T_{n}(i)}{V_{n}(i)}\cdot I_{\{V_{n}(i)\geq 1\}}+
w\cdot I_{\{V_{n}(i)=0\}}.
\end{eqnarray*}
We also define 
\begin{eqnarray}\label{defloccluster}
C(i)
=
\frac{ E_{T}(X_{i}) }
{ E_{D}(X_{i})^{2} }
\cdot I_{\{E_{D}(X_{i})>0\}}+
w\cdot I_{\{E_{D}(X_{i})=0\}}, 
\end{eqnarray}
where 
\begin{eqnarray*}
& &E_{D}(X_{i})=\int _{\mathbb{R}} h_{D}(X_{i},y)F(dy),\\
& &E_{T}(X_{i})=\int _{\mathbb{R}}\int _{\mathbb{R}} h_{T}(X_{i},y,z)F(dy)F(dz).
\end{eqnarray*}
We consider 
$C_{n}(i;x)$, $C(i;x)$, $D_{n}(i;x)$, $T_{n}(i;x)$, and $V_{n}(i;x)$, which are
random variables $C_{n}(i)$, $C(i)$, $D_{n}(i)$, $T_{n}(i)$ and $V_{n}(i)$
restricted to the subspace such that $\{X_{i}=x\}$. 
For example, $T_{n}(1;x)=\sum _{2\leq j<k\leq n}h_{T}(x,X_{j},X_{k})$. 
We obtain the following asymptotic results for $C_{n}(i)$: 
\begin{theorem}\label{thmloccluster}
As $n\to \infty $, 
\newline 
(i)\ For any $x\in \mathbb{R}$, $C_{n}(1;x)\to C(1;x)$, almost surely. 
\newline 
In particular, 
\newline 
(ii)\ $C_{n}(1)\to C(1)$, almost surely.
\end{theorem}
\begin{proof}
For an arbitrary fixed $x\in \mathbb{R}$, we first prove
\begin{eqnarray}\label{iffzero}
\mathbb{E}\left[h_{V}(x,X_{2},X_{3})\right]=0 
\Longleftrightarrow 
\mathbb{P}\left(V_{n}(1;x)=0 \text{ for all }n\geq 1\right)=1.
\end{eqnarray}
Indeed, if $\mathbb{E}\left[h_{V}(x,X_{2},X_{3})\right]=0$ then 
\begin{eqnarray*}
\mathbb{E}\left[V_{n}(1;x)\right]
=
\mathbb{E}\left[\sum _{2\leq j<k\leq n}h_{V}(x,X_{j},X_{k})\right]
=
\sum _{2\leq j<k\leq n}\mathbb{E}\left[h_{V}(x,X_{j},X_{k})\right]
=0
\end{eqnarray*}
for all $n\geq 1$. Conversely, if $\mathbb{E}\left[V_{n}(1;x)\right]=0$ 
for all $n\geq 1$ then 
\begin{eqnarray*}
\mathbb{E}\left[h_{V}(x,X_{2},X_{3})\right]
=
\mathbb{E}\left[V_{3}(1;x)\right]
=0.
\end{eqnarray*}
Therefore, we obtain
\begin{eqnarray*}
\mathbb{E}\left[h_{V}(x,X_{2},X_{3})\right] =0
\Longleftrightarrow 
\mathbb{E}\left[V_{n}(1;x)\right]=0
\text{ for all } n\geq 1.
\end{eqnarray*}
Since $V_{n}(1;x)$ is nonnegative, 
\begin{eqnarray*}
\mathbb{E}\left[V_{n}(1;x)\right]=0
\Longleftrightarrow 
\mathbb{P}\left(V_{n}(1;x)=0\right)=1
\text{ for all } n\geq 1. 
\end{eqnarray*}
Moreover, $\{V_{n}(1;x)=0\}$ is nonincreasing with $n$, which implies  
\begin{eqnarray*}
\mathbb{P}\left(V_{n}(1;x)=0\right)=1 \text{ for all }n\geq 1
\Longleftrightarrow 
\mathbb{P}\left(V_{n}(1;x)=0\text{ for all }n\geq 1\right)=1.
\end{eqnarray*}
Thus we have Eq.\ (\ref{iffzero}). 
By definition, $V_{n}(1;x)$ is invariant under any permutation on 
$\{x_{2},x_{3},\ldots ,x_{n}\}$, 
and $V_{n}(1;x)$ is nondecreasing, i.e., 
\begin{eqnarray}\label{nondecrease}
V_{n}(1;x)(x_{2},x_{3},\ldots ,x_{n})\leq 
V_{n+1}(1;x)(x_{2},x_{3},\ldots ,x_{n},x_{n+1}) 
\end{eqnarray}
for all $n\geq 1$.  
Therefore $\mathbb{P}\left(V_{n}(1;x)=0 \text{ for all }n\geq 1\right)$
equals to zero or one by the Hewitt-Savage zero-one law (see Theorem 36.5 of \cite{b}). 
So we have 
\begin{eqnarray*}
\mathbb{E}\left[h_{V}(x,X_{2},X_{3})\right]>0 
\Longleftrightarrow 
\mathbb{P}\left(V_{n}(1;x)=0 \text{ for all }n\geq 1\right)=0
\Longleftrightarrow 
\mathbb{P}\left(V_{n}(1;x)\geq 1 \text{ for some }n\geq 1\right)=1.
\end{eqnarray*}
Using Eq. (\ref{nondecrease}), 
$\left\{V_{n}(1;x)\geq 1 \text{ for some }n\geq 1\right\}$ 
is equivalent to the event \\
$\left\{\exists N\geq 1 \text{ s.t. }V_{n}(1;x)\geq 1 \text{ for all }n\geq N\right\}$. 
Hence we obtain
\begin{eqnarray}\label{iffone}
\mathbb{E}\left[h_{V}(x,X_{2},X_{3})\right]>0 
\! \! \! \! &\Longleftrightarrow&\! \! \! \! 
\mathbb{P}\left(\exists N\geq 1 \text{ s.t. }V_{n}(1;x)\geq 1 \text{ for all }n\geq N\right)=1\notag \\
\! \! \! \! &\Longleftrightarrow&\! \! \! \! 
\mathbb{P}\left(\exists N\geq 1 \text{ s.t. }C_{n}(1;x)=T_{n}(1;x)/V_{n}(1;x) \text{ for all }n\geq N\right)=1.
\end{eqnarray}
Since $h_{T}(x,x_{2},x_{3})$ and $h_{V}(x,x_{2},x_{3})$ are 
symmetric functions of $x_{2}$ and $x_{3}$, we define $U$-statistics 
\begin{eqnarray*}
\frac{ T_{n}(1;x) }{ \binom{n-1}{2} }
\! \! \! \! &=&\! \! \! \! 
\frac{1}{ \binom{n-1}{2} }\sum _{2\leq j<k\leq n}h_{T}(x,X_{j},X_{k}),\\
\frac{ V_{n}(1;x) }{ \binom{n-1}{2} }
\! \! \! \! &=&\! \! \! \! 
\frac{1}{ \binom{n-1}{2} }\sum _{2\leq j<k\leq n}h_{V}(x,X_{j},X_{k}).
\end{eqnarray*}
We have the following SLLN by Theorem A in Section 5.4 of \cite{s}: As $n\to \infty $,
\begin{eqnarray}
\frac{ T_{n}(1;x) }{ \binom{n-1}{2} }
\! \! \! \! &\to &\! \! \! \! 
\mathbb{E}\left[h_{T}(x,X_{2},X_{3})\right],
\quad \text{almost surely,}\label{asT}\\
\frac{ V_{n}(1;x) }{ \binom{n-1}{2} }
\! \! \! \! &\to &\! \! \! \! 
\mathbb{E}\left[h_{V}(x,X_{2},X_{3})\right],
\quad \text{almost surely.}\label{asV}
\end{eqnarray}
Based on Eqs.\ (\ref{asT}) and (\ref{asV}),
the corresponding clustering coefficient 
\begin{eqnarray*}
C_{n}(1;x)=
\frac{T_{n}(1;x)}{V_{n}(1;x)}=\frac{T_{n}(1;x)/\binom{n-1}{2}}{V_{n}(1;x)/\binom{n-1}{2}}
\end{eqnarray*}
converges to $\mathbb{E}\left[h_{T}(x,X_{2},X_{3})\right]/\mathbb{E}\left[h_{V}(x,X_{2},X_{3})\right]$, 
almost surely as $n\to \infty $. By Eq.\ (\ref{iffone}), we have 
\begin{eqnarray}\label{convone}
\mathbb{P}\left(\lim _{n\to \infty }C_{n}(1;x)=
\mathbb{E}\left[h_{T}(x,X_{2},X_{3})\right]/\mathbb{E}\left[h_{V}(x,X_{2},X_{3})\right]\right)=1.
\end{eqnarray}
On the other hand, Eqs.\ (\ref{iffzero}) and (\ref{convone}) 
imply that 
$\mathbb{E}\left[h_{V}(x,X_{2},X_{3})\right]=0$ is equivalent to \\
$\mathbb{P}\left(\lim _{n\to \infty }C_{n}(1;x)=w\right)=1$.
With the relation $\mathbb{E}\left[h_{V}(x,X_{2},X_{3})\right]=\mathbb{E}\left[h_{D}(x,X_{2})\right]^{2}$, 
we obtain
\begin{eqnarray*}
C_{n}(1;x)
\to 
C(1;x)
=
\frac{\mathbb{E}\left[h_{T}(x,X_{2},X_{3})\right]}{\mathbb{E}\left[h_{D}(x,X_{2})\right]^{2}}
\cdot I_{\left\{\mathbb{E}\left[h_{D}(x,X_{2})\right]>0\right\}}
+
w
\cdot I_{\left\{\mathbb{E}\left[h_{D}(x,X_{2})\right]=0\right\}},
\end{eqnarray*}
almost surely as $n\to \infty $.
Particularly, we have by using Fubini's theorem, 
\begin{eqnarray*}
\mathbb{P}\left(\lim _{n\to \infty }C_{n}(1)=C(1)\right)
=
\int _{\mathbb{R}}
\mathbb{P}\left(\lim _{n\to \infty }C_{n}(1;x)=C(1;x)\right)F(dx)
=
\int _{\mathbb{R}}1\cdot F(dx)=1.
\end{eqnarray*} 
This completes the proof. 
\end{proof}

\subsection{Global Clustering Coefficient}

The global clustering coefficient is defined by 
\begin{eqnarray*}
C_{n}=\frac{1}{n}\sum _{i=1}^{n}C_{n}(i). 
\end{eqnarray*}
Since it is a symmetric function of $(x_{1},\dots ,x_{n})$, 
we can prove SLLN for $C_n$ by using the ergodic theory. 
\begin{theorem}\label{thmglobcluster}
As $n\to \infty $, 
\begin{eqnarray*}
C_{n}\to \mathbb{E}\left[C(1)\right],
\quad \text{almost surely}.
\end{eqnarray*} 
\end{theorem}
\begin{proof}
For simplicity, we only deal with the case $\mathbb{E}\left[C(1)\right]=0$. 
For general cases, we can prove the theorem by replacing $C(1)$ by $C(1)-\mathbb{E}\left[C(1)\right]$. 
Let $\boldsymbol{x}=(x_{1},x_{2},\ldots)$ be an infinite vector and $\boldsymbol{x}_{k}=x_{k}$. 
We define measure-preserving transformation $T_{n}$ for the product measure $\mathbb{P}$ by 
\begin{eqnarray*}
\left(T_{n}\boldsymbol{x}\right)_{k}=
\begin{cases}
x_{k+1}& \text{if $1\leq k\leq n-1$},\\
x_{1}& \text{if $k=n$},\\
x_{k}& \text{otherwise},
\end{cases}
\end{eqnarray*}
for each $n\geq 1$. By denoting $C_{n}(i;\boldsymbol{x})=C_{n}(i;x_{i})$,
a realization 
of $C_{n}$ is represented by 
\begin{eqnarray*}
C_{n}(\boldsymbol{x})=\frac{1}{n}\sum _{i=0}^{n-1}C_{n}(1;T_{n}^{i}\boldsymbol{x}).
\end{eqnarray*}
For arbitrary fixed $\varepsilon >0$, we define 
\begin{eqnarray}\label{eqstar}
C_{n}^{\varepsilon }(1;\boldsymbol{x})
\! \! \! \! &=&\! \! \! \! 
\left(C_{n}(1;\boldsymbol{x})-\varepsilon \right)\cdot I_{A_{\varepsilon }}, \\
\ S_{n}^{\varepsilon }(\boldsymbol{x})
\! \! \! \! &=&\! \! \! \! 
\sum _{i=0}^{n-1}C_{n}^{\varepsilon }(1;T_{n}^{i}\boldsymbol{x}), \notag 
\end{eqnarray}
where 
$
A_{\varepsilon }=\left\{\boldsymbol{x}:\limsup _{n\to \infty }C_{n}(\boldsymbol{x})>\varepsilon \right\}.
$
Using the maximal ergodic theorem (see Theorem 24.2 of \cite{b}), 
\begin{eqnarray*}
\int _{M_{n}^{\varepsilon}}C_{n}^{\varepsilon }(1;\boldsymbol{x})d\mathbb{P}\geq 0
\end{eqnarray*}
for every $n\geq 1$, 
where $M_{n}^{\varepsilon }=\left\{\boldsymbol{x}:\sup _{1\leq j\leq n}S_{j}^{\varepsilon }(\boldsymbol{x})>0\right\}$. 
On the other hand, we have 
\begin{eqnarray*}
M_{n}^{\varepsilon }\uparrow 
\left\{\boldsymbol{x}:\sup _{k\geq 1}S_{k}^{\varepsilon }(\boldsymbol{x})>0\right\}
=
\left\{\boldsymbol{x}:\sup _{k\geq 1}\frac{ S_{k}^{\varepsilon }(\boldsymbol{x}) }{k}>0\right\}
=
\left\{\boldsymbol{x}:\sup _{k\geq 1}C_{n}(\boldsymbol{x})>\varepsilon \right\}\cap A_{\varepsilon }
=
A_{\varepsilon },
\end{eqnarray*}
as $n\to \infty $ by Eq.\ (\ref{eqstar}). 
From the
dominated convergence theorem and Theorem \ref{thmloccluster}, we derive
\begin{eqnarray}\label{eqlimint}
0\leq \int _{M_{n}^{\varepsilon }}C_{n}^{\varepsilon }(1;\boldsymbol{x})d\mathbb{P}
\to 
\int _{A_{\varepsilon }}\left[C(1;\boldsymbol{x})-\varepsilon \right]d\mathbb{P},
\end{eqnarray}
as $n\to \infty $. 

Let $\mathcal{I}_{n}$ be the class of sets that are invariant under all permutations of 
the first $n$ coodinates and $\mathcal{I}=\bigcap _{n=1}^{\infty }\mathcal{I}_{n}$. 
It is easy to check that $A_{\varepsilon }\in \mathcal{I}$. 
Since $\mathbb{P}(A)$ equals to zero or one for any $A\in \mathcal{I}$ 
by the Hewitt-Savage zero-one law (see Theorem 36.5 of \cite{b}), 
the conditional expectation $\mathbb{E}\left[C(1)|\mathcal{I}\right]$ equals to 
$\mathbb{E}\left[C(1)\right]$, almost surely. This leads to 
\begin{eqnarray*}
0\leq \int _{A_{\varepsilon }}\left[C(1;\boldsymbol{x})-\varepsilon \right]d\mathbb{P}
\! \! \! \! &=&\! \! \! \! 
\int _{A_{\varepsilon }}C(1;\boldsymbol{x})d\mathbb{P}
-\varepsilon\mathbb{P}(A_{\varepsilon })
=
\int _{A_{\varepsilon }}\mathbb{E}\left[C(1;\boldsymbol{x})|\mathcal{I}\right]d\mathbb{P}
-\varepsilon\mathbb{P}(A_{\varepsilon })\\
\! \! \! \! &=&\! \! \! \! 
\int _{A_{\varepsilon }}\mathbb{E}\left[C(1)\right]d\mathbb{P}
-\varepsilon\mathbb{P}(A_{\varepsilon })
=
-\varepsilon\mathbb{P}(A_{\varepsilon }), 
\end{eqnarray*}
by Eq.\ (\ref{eqlimint}) and $\mathbb{E}\left[C(1)\right]=0$. 
Then we have $\mathbb{P}(A_{\varepsilon })=0$ for any $\varepsilon >0$ and therefore 
$\limsup _{n\to \infty }C_{n}\leq 0$, almost surely. 
Repeating the same argument for $-C_{n}$, 
we have $\liminf _{n\to \infty }C_{n}\geq 0$, almost surely. 
This completes the proof. 
\end{proof}
Here we show a simple example for Theorem \ref{thmglobcluster}.
Consider an i.i.d. sequence $X_{1}, \ldots ,X_{n}$ such that $\mathbb{P}(X_{i}=1)=p$ and 
$\mathbb{P}(X_{i}=0)=1-p$ for all $i=1,\ldots ,n$. 
Let $B_{1}=(\theta ,\infty )$ and $f_{c}(x,y)=x+y$. 
We set a threshold $\theta $ such that $0<\theta <1$. 
In this case, a pair of vertices $i$ and $j$ with $i\neq j$ is disconnected 
if and only if $X_{i}=X_{j}=0$. By direct computation, we have 
\begin{eqnarray*}
\mathbb{E}\left[C(1)\right]
=
p\cdot C(1;1)+(1-p)\cdot C(1;0)
=p\cdot \frac{p^{2}+2p(1-p)}{1}+(1-p)\cdot \frac{p^{2}}{p^{2}}
=
1-p(1-p)^{2}. 
\end{eqnarray*}
In order to calculate $C_{n}$, let $S_{n}=\sum _{i=1}^{n}X_{i}$, that
is, the number of vertices with $X_{i}=1$. 
We use the symbols $x_{i}$, $s_{n}$ and $c_{n}$ as realization of random variables 
$X_{i}$, $S_{n}$ and $C_{n}$ respectively. 
If $s_{n}=0$, the graph consists of $n$ isolated vertices. 
In this case $c_{n}=w$. 
If $s_{n}=1$, the graph is the star in which 
only one central vertex has $n-1$ edges and other $n-1$ vertices are
connected only to the center. So we obtain
\begin{eqnarray*}
c_{n}=\frac{1}{n}\left\{0\cdot 1+w\cdot (n-1)\right\}=\left(1-\frac{1}{n}\right)\cdot w.
\end{eqnarray*}
If $2\leq s_{n}\leq n-2$, $s_{n}$ vertices with $x_{i}=1$ have
 $n-1$ edges, 
and the other $n-s_{n}$ vertices are connected only to
the vertices with $x_{i}=1$. So we have
\begin{eqnarray*}
c_{n}
=
\frac{1}{n}
\left\{
\frac{ \binom{n-1}{2}-\binom{n-s_{n}}{2} }{ \binom{n-1}{2} }\cdot s_{n}
+
\frac{ \binom{s_{n}}{2} }{ \binom{s_{n}}{2} }\cdot (n-s_{n})
\right\}
= 
1-\frac{(n-s_{n})(n-1-s_{n})s_{n}}{n(n-1)(n-2)}.
\end{eqnarray*}
If $s_{n}=n-1$ or $n$, we obtain the complete graph, and $c_{n}=1$. 
Noting
\begin{eqnarray*}
1-\frac{(n-S_{n})(n-1-S_{n})S_{n}}{n(n-1)(n-2)}
=
\begin{cases}
1&\text{if}\ \ s_{n}=0,n-1,n,\\
1-\frac{1}{n}&\text{if}\ \ s_{n}=1,\\
1-\frac{(n-s_{n})(n-1-s_{n})s_{n}}{n(n-1)(n-2)}&\text{otherwise},
\end{cases}
\end{eqnarray*}
we have 
\begin{eqnarray*}
C_{n}
\! \! \! \! &=&\! \! \! \! 
\left[1-\frac{(n-S_{n})(n-1-S_{n})S_{n}}{n(n-1)(n-2)}\right]\cdot 
\left\{I_{\{2,\ldots ,n\}}(S_{n})+r\cdot I_{\{0,1\}}(S_{n})\right\}\\
\! \! \! \! &=&\! \! \! \!  
\left[1-\left(1-\frac{ S_{n} }{n}\right)\left(1-\frac{ S_{n} }{n-1}\right)
\left(\frac{ S_{n} }{n-2}\right)\right]\cdot 
\left\{1+(r-1)\cdot I_{\{0,1\}}(S_{n})\right\}\\
\! \! \! \! &\to &\! \! \! \!  
1-(1-p)^{2}p=\mathbb{E}\left[C(1)\right],\quad \text{almost surely}\ (n\to \infty ).
\end{eqnarray*}
The last convergence comes from SLLN for the i.i.d. sequence. 

One of our motivations to study limit theorems for the clustering 
coefficients is to make a clear distinction between the proportion of 
triangles in an entire graph and the clustering coefficients.
 By Eq.\ (66) of \cite{kmrs05}, the normalized number of triangles 
including vertex $1$ converges to $E_{T}(x)$, almost surely for each 
realization $x$ of $X_{1}$, where the normalization constant is 
equal to $\binom {n-1}{2}$. 
In the same way, the degree of vertex $1$ normalized by $n-1$ 
converges to $E_{D}(x)$, almost surely. Thus, when $E_{D}(x)>0$, the 
local clustering coefficient converges almost surely to 
$E_{T}(x)/E_{D}(x)^{2}$, that is, the limit of the normalized number of 
triangles divided by the square of the limit of the normalized degree. 
The mean field result corresponding to Theorem \ref{thmloccluster} is 
found in Eq.\ (3) of \cite{scb}. The denominator equals to $E_{T}(x)$ 
and the numerator $[k(x)/N]^{2}$ converges to $E_{D}(x)^{2}$, almost 
surely as $N\to \infty $, where $k(x)$ is the degree of the vertex $1$ 
and $N$ denotes the number of vertices. 
Equation (30) of \cite{bps} corresponds to the 
normalized number of triangles. These 
heuristic results are consistent with our rigorous result. 
Several examples for the global clustering coefficient are calculated 
in \cite{mmk04}.

In practice, we may substitute $0$ into $w$ and consider 
\begin{eqnarray*}
\Tilde{C}_{n}
\! \! \! \! &=&\! \! \! \! 
\frac{1}{n-\text{number of vertices with degree $0$ or $1$}}\ \sum _{i=1}^{n}C_{n}(i)\\
\! \! \! \! &=&\! \! \! \! 
\frac{1}{n-\sum _{i=1}^{n}I_{\{0\}}(V_{n}(i))}\ \sum _{i=1}^{n}C_{n}(i), 
\end{eqnarray*}
instead of $C_{n}$. 
Using the same arguments of 
Theorems \ref{thmloccluster} and \ref{thmglobcluster}, 
it is easy to prove that
\begin{eqnarray*}
\lim _{n\to \infty }\frac{1}{n}\sum _{i=1}^{n}I_{\{0\}}(V_{n}(i))
=
\lim _{n\to \infty }\frac{1}{n}\sum _{i=1}^{n}I_{\{w\}}(C_{n}(i))
=
\mathbb{P}(C(1)=w)
=
\mathbb{P}(E_{D}(X_{1})=0),\quad \text{almost surely.}
\end{eqnarray*}
The last equality follows from the
definition of $C(1)$, i.e., Eq.\ (\ref{defloccluster}).
Noting that 
\begin{eqnarray*}
\Tilde{C}_{n}=\frac{1}{1-(1/n)\sum _{i=1}^{n}I_{\{0\}}(V_{n}(i))}\cdot C_{n}, 
\end{eqnarray*}
we have the following:
\begin{corollary}
As $n\to \infty $, 
\begin{eqnarray*}
\Tilde{C}_{n}\to \frac{1}{1-\mathbb{P}(E_{D}(X_{1})=0)}\cdot \mathbb{E}\left[C(1)\right],
\quad \text{almost surely}.
\end{eqnarray*}
\end{corollary}
\section{Examples}
In this section, we show examples of the 
limit degree distribution, 
i.e., $m=2$ and $\mathcal{A}_{2}$ is chosen as 
the collection of all possible edges
in the limit theorem (Fact \ref{thmlocu}).
We consider the case $l=1$ with $f_c\equiv f_c^1$.
We assume that the random variable $X_{1}$ 
is absolutely continuous so that it 
has a probability density function $f$.
Let $\supp f=\overline{\{x\in \mathbb{R} : f(x)\neq 0\}}$ be the support of $f$.

We first set $B_1=(\theta ,\infty )$ and $f_{c}(x,y)=x+y$, i.e.,
the threshold network model in which
an edge $\langle i,j\rangle $ forms if 
$\theta <X_{i}+X_{j}$ for a given threshold $\theta \in \mathbb{R}$ 
\cite{bps,ccdm,kmrs05,mmk04}.
By calculating the characteristic function of 
$D=U(B_1,\mathcal{A}_{2})$, namely, the
density of edges connected to vertex $1$, 
we obtain the following results: 
\begin{lemma}
\begin{enumerate}
\item
If there exists $a\in \mathbb{R}$ such that 
$\supp f=[a,\infty )$, then 
\begin{eqnarray*}
D\thicksim 
\begin{cases}
\delta _{1}(dk)& \text{{\rm if}\quad $\theta \leq 2a$,}\\
I_{(1-F(\theta -a),1)}(k)\cdot \frac{f\bigl(\theta -F^{-1}(1-k)\bigr)}{f\bigl(F^{-1}(1-k)\bigr)}dx\\
+\bigl(1-F(\theta -a)\bigr)\cdot \delta _{1}(dk)
& \text{{\rm if}\quad $\theta >2a$.}
\end{cases}
\end{eqnarray*}
\item
If there exists $b\in \mathbb{R}$ such that 
$\supp f=(-\infty ,b]$, then 
\begin{eqnarray*}
D\thicksim 
\begin{cases}
\bigl(1-F(\theta -b)\bigr)\cdot \delta _{0}(dk)\\
+I_{(0,1-F(\theta -b))}(k)\cdot \frac{f\bigl(\theta -F^{-1}(1-k)\bigr)}{f\bigl(F^{-1}(1-k)\bigr)}dk
& \text{{\rm if}\quad $\theta <2b$,}\\
\delta _{0}(dk)& \text{{\rm if}\quad $\theta \geq 2b$.}
\end{cases}
\end{eqnarray*}
\item
If there exist $a,b\in \mathbb{R}$ such that 
$\supp f=[a,b]$, then 
\begin{eqnarray*}
D\thicksim 
\begin{cases}
\delta _{1}(dk)& \text{{\rm if}\quad $\theta \leq 2a$,}\\
I_{(1-F(\theta -a),1)}(k)\cdot \frac{f\bigl(\theta -F^{-1}(1-k)\bigr)}{f\bigl(F^{-1}(1-k)\bigr)}dk\\
+\bigl(1-F(\theta -a)\bigr)\cdot \delta _{1}(dk)
& \text{{\rm if}\quad $2a<\theta <a+b$,}\\
I_{(0,1)}(k)\cdot \frac{f\bigl(a+b-F^{-1}(1-k)\bigr)}{f\bigl(F^{-1}(1-k)\bigr)}dk
& \text{{\rm if}\quad $\theta =a+b$,}\\
\bigl(1-F(\theta -b)\bigr)\cdot \delta _{0}(dk)\\
+I_{(0,1-F(\theta -b))}(k)\cdot \frac{f\bigl(\theta -F^{-1}(1-k)\bigr)}{f\bigl(F^{-1}(1-k)\bigr)}dk
& \text{{\rm if}\quad $a+b<\theta <2b$,}\\
\delta _{0}(dk)& \text{{\rm if}\quad $\theta \geq 2b$.}
\end{cases}
\end{eqnarray*}
Furthermore, if $f$ is symmetric on $\supp f$, then 
\begin{eqnarray*}
D\thicksim I_{(0,1)}(k)dk\quad \text{{\rm if}\quad $\theta =a+b$.}
\end{eqnarray*}
\item
If $\supp f=(-\infty ,\infty )$, then 
\begin{eqnarray*}
D\thicksim 
I_{(0,1)}(k)\cdot \frac{f\bigl(\theta -F^{-1}(1-k)\bigr)}{f\bigl(F^{-1}(1-k)\bigr)}dk
\end{eqnarray*}
for any $\theta \in \mathbb{R}$.
\end{enumerate}
\end{lemma}

\begin{example}({\it Exponential distribution})
If the random variable $X_{1}$ has the probability density function 
\begin{eqnarray}\label{densexp}
f(x)=
\begin{cases}
\lambda e^{-\lambda x}& \text{if\quad $x\geq 0$,}\\
0& \text{otherwise,}
\end{cases}
\end{eqnarray}
for a given $\lambda >0$, then 
\begin{eqnarray*}\label{lim1}
D\thicksim 
\begin{cases}
\delta _{1}(dk)& \text{if\quad $\theta \leq 0$,}\\
I_{(e^{-\lambda \theta },1)}(k)\cdot \frac{e^{-\lambda \theta }}{k^{2}}dk
+e^{-\lambda \theta }\cdot \delta _{1}(dk)
& \text{if\quad $\theta >0$.}
\end{cases}
\end{eqnarray*}
\end{example}
\begin{example}({\it Pareto distribution})
If  
\begin{eqnarray*}
f(x)=
\begin{cases}
\frac{c}{a}\cdot \bigl(\frac{a}{x}\bigr)^{c+1}& \text{if\quad $x\geq a$,}\\
0& \text{otherwise,}
\end{cases}
\end{eqnarray*}
for given $a,c>0$, then 
\begin{eqnarray*}
D\thicksim 
\begin{cases}
\delta _{1}(dk)& \text{if\quad $\theta \leq 2a$,}\\
I_{((\frac{a}{\theta -a})^{c},1)}(k)\cdot \bigl(\frac{a}{\theta \cdot k^{1/c}-a}\bigr)^{c+1}dk
+\bigl(\frac{a}{\theta -a}\bigr)^{c}\cdot \delta _{1}(dk)
& \text{if\quad $\theta >2a$.}
\end{cases}
\end{eqnarray*}
\end{example}
The distribution of $D$ of these two examples is 
proportional to $k^{-\alpha}$. 
The exponent $\alpha $ equals $2$ in 
Example 1 and $1+1/c$ in Example 2.
Because of a lower cutoff of $f$ in both examples,
the limit distributions have weights on $\delta _{1}$. 
\begin{example}({\it Uniform distribution})
If 
\begin{eqnarray*}
f(x)=
\begin{cases}
1& \text{if\quad $0\leq x\leq 1$,}\\
0& \text{otherwise,}
\end{cases}
\end{eqnarray*}
then 
\begin{eqnarray*}
D\thicksim 
\begin{cases}
\delta _{1}(dk)& \text{if\quad $\theta \leq 0$,}\\
I_{(1-\theta ,1)}(k)dk
+(1-\theta )\cdot \delta _{1}(dk)
& \text{if\quad $0<\theta <1$,}\\
I_{(0,1)}(k)dk
& \text{if\quad $\theta =1$,}\\
(\theta -1)\cdot \delta _{0}(dk)
+I_{(0,2-\theta )}(k)dk
& \text{if\quad $1<\theta <2$,}\\
\delta _{0}(dk)& \text{if\quad $\theta \geq 2$.}
\end{cases}
\end{eqnarray*}
\end{example}
In this case, the limit distribution is mixture of the
uniform distribution and the delta measure. 

By choosing $B_1=(\theta _{1},\theta _{2}]$ and $f_{c}(x,y)=x+y$, 
we obtain a
generalization of the model investigated in \cite{bps,ccdm,kmrs05,mmk04}. 
More precisely,
an edge $\langle i,j\rangle $ forms
if $\theta _{1}<X_{i}+X_{j}\leq \theta _{2}$ for given 
thresholds $\theta _{1},\theta _{2}\in \mathbb{R}$ such that $\theta _{1}<\theta _{2}$. 
To calculate the characteristic function of 
$D=U(B_1, \mathcal{A}_{2})$, 
we consider the case in which
a random variable $X_{1}$ has the probability density function (\ref{densexp}), 
i.e.,\ the exponential distribution, for which
the limit distribution is represented by: 
\begin{eqnarray*}
D&\thicksim &
\begin{cases}
\delta _{0}(dk)
& \text{if\quad $\theta _{1}<\theta _{2}\leq 0$,}\\
e^{-\lambda \theta _{2}}\cdot \delta _{0}(dk)\\
+I_{(0,1-e^{-\lambda \theta _{2}})}(k)\cdot \frac{e^{-\lambda \theta _{2}}}{(1-k)^{2}}dk 
& \text{if\quad $\theta _{1}\leq 0<\theta _{2}$,}\\
e^{-\lambda \theta _{2}}\cdot \delta _{0}(dk)\\
+I_{(0,1-e^{-\lambda (\theta _{2}-\theta _{1})})}(k)\cdot \frac{e^{-\lambda \theta _{2}}}{(1-k)^{2}}dk\\
+I_{(e^{-\lambda \theta _{1}}-e^{-\lambda \theta _{2}},1-e^{-\lambda (\theta _{2}-\theta _{1})})}(k)\cdot \\
\frac{e^{-\lambda \theta _{1}}-e^{-\lambda \theta _{2}}}{k^{2}}dk
& \text{if\quad $0<\theta _{1}<\theta _{2}$.}
\end{cases}
\end{eqnarray*}

Finally, we deal with an example with $l=2$. 
For fixed $\theta \in \mathbb{R}$ and $c\in [0,\infty )$, 
we choose $B_{1}=(\theta ,\infty ], B_{2}=(0,c]$,
$f_{c}^{1}(x,y)=x+y$, and $f_{c}^{2}(x,y)=|x-y|$. 
We consider the case in which $X_{1}$ is distributed according
to the exponential distribution (Eq.\ (\ref{densexp})). This is the model proposed in \cite{mk06}. 
Because the kernel function of this model is  
\begin{eqnarray*}
h(x,x_{2})=
\begin{cases}
I_{[-c+x,c+x]}(x_{2})
& \text{if\quad $\frac{\theta +c}{2}\leq x$,}\\
I_{(\theta -x,c+x]}(x_{2})
& \text{if\quad $\frac{\theta -c}{2}\leq x\leq \frac{\theta +c}{2}$,}\\
0
& \text{if\quad $x\leq \frac{\theta -c}{2}$,}\\
\end{cases}
\end{eqnarray*}
the limit distribution $D=U(\mathcal{C}_{\theta ,c},\mathcal{A}_{2})$ is the following: 
\begin{eqnarray*}
D&\thicksim &
\begin{cases}
\left(1-e^{-\lambda (\theta -c)/2}\right)\delta _{0}(dk)\\
+
I_{(0,2e^{-\lambda (\theta +c)/2}\sinh (\lambda c)]}(k)\cdot g(k)dk
& \text{if\quad $c\leq \theta $,}\\
I_{(0,e^{-\lambda \theta }-e^{-\lambda c}]}(k)\cdot \frac{1}{2\sinh (\lambda c)}dk\\
+
I_{(e^{-\lambda \theta }-e^{-\lambda c},1-e^{-\lambda (\theta +c)}]}(k)\cdot g(k)dk\\
+
I_{(1-e^{-\lambda (\theta +c)},1-e^{-2\lambda c}]}(k)\cdot \frac{e^{2\lambda c}}{2\sinh (\lambda c)}dk
& \text{if\quad $0\leq \theta \leq c$,}\\
I_{(0,1-e^{-\lambda c}]}(k)\cdot \frac{1}{2\sinh (\lambda c)}dk\\
+
I_{(1-e^{-\lambda c},1-e^{-2\lambda c}]}(k)\cdot \frac{e^{2\lambda c}}{2\sinh (\lambda c)}dk
& \text{if\quad $-c\leq \theta \leq 0$,}\\
I_{(0,1-e^{-\lambda c}]}(k)\cdot \frac{1}{2\sinh (\lambda c)}dk\\
+
I_{(1-e^{-\lambda c},1-e^{-2\lambda c}]}(k)\cdot \frac{e^{2\lambda c}}{2\sinh (\lambda c)}dk
& \text{if\quad $\theta \leq -c$,}
\end{cases}
\end{eqnarray*}
where 
\begin{eqnarray*}
g(k)=
\frac{ 4e^{-\lambda \theta } }
{ \left(k+\sqrt{k^{2}+4e^{-\lambda (\theta +c)}}\right)^{2}+4e^{-\lambda (\theta +c)} }
+
\frac{1}{2\sinh (\lambda c)}.
\end{eqnarray*}
\begin{acknowledgements}
{\rm We thank Masato Takei for valuable discussions and comments.}
\end{acknowledgements}

\end{document}